# REMARKS ON THE $\mathcal{O}(N)$ IMPLEMENTATION OF THE FAST MARCHING METHOD

CHRISTIAN RASCH AND THOMAS SATZGER

ABSTRACT. The fast marching algorithm computes an approximate solution to the eikonal equation in $\mathcal{O}(N \log N)$ time, where the factor $\log N$ is due to the administration of a priority queue. Recently, [YBS06] have suggested to use an untidy priority queue, reducing the overall complexity to $\mathcal{O}(N)$ at the price of a small error in the computed solution. In this paper, we give an explicit estimate of the error introduced, which is based on a discrete comparison principle. This estimates implies in particular that the choice of an accuracy level that is independent of the speed function $F$ results in the complexity bound $\mathcal{O}(F_{\max}/F_{\min} N)$. A numerical experiment illustrates this robustness problem for large ratios $F_{\max}/F_{\min}$.

## 1. INTRODUCTION

The fast marching method, introduced in [Set96], is an efficient algorithm to compute the discrete solution to the eikonal equation $\|\nabla T(x)\| F(x) = 1$, with speed function $F$. The algorithm is closely related to Dijkstra's single-source shortest path algorithm, and computes the grid-solution, starting with the grid-points adjacent to the boundary, by traversing the computational domain along increasing values of $T$. The total complexity of this method is $\mathcal{O}(N \log N)$, where $N$ denotes the number of grid-points. Here, the factor $\log N$ comes from the administration of a priority queue.

In [Tsi95], Tsitsiklis suggested independently a similar algorithm. Other than [Set96], the discretization of the eikonal equation in [Tsi95] was obtained in the framework of optimal control theory. Moreover, Tsitsiklis showed that a bucket sort technique, together with a slightly different discretization, may be used to compute an approximate solution to the eikonal equation in $\mathcal{O}(N)$ time.

Another approach has recently been suggested by [YBS06]: using an *untidy* priority queue in the fast marching method also reduces the complexity to optimal $\mathcal{O}(N)$. The idea is to use a bucket sort technique together with a quantization that does not distinguish between values of $T$ within a small range. This way, while sorting the values of $T$ becomes less expensive, an error is introduced, however, which should be of the same order as the local truncation error of the discretization. The central result of this paper is Theorem 4, which gives a rigorous estimate of the error introduced. Before we come to the theorem, the discretization and some basic notation is introduced in Section 2. Then, in Section 3, we briefly recall the fast marching method and the bucket sort technique, as proposed in [YBS06]. A discrete comparison principle and the error estimation are subject of Section 4. Finally, an example is given, which illustrates the theoretical result.

## 2. Discretization

Given a closed subset $\Gamma \subset [0,1]^2$, we consider the eikonal equation

$$\|\nabla T(x)\| \, F(x) = 1 \quad (x \in [0,1]^2 \setminus \Gamma), \qquad T|_\Gamma = 0, \tag{1}$$

with continuous speed function $F : [0,1]^2 \to \mathbb{R}$, such that there are real numbers $F_{\min}, F^{\max}$ with

$$0 < F_{\min} \leqslant F(x) \leqslant F_{\max} \qquad \forall x. \tag{2}$$

The computational domain $[0,1]^2$ is endowed with a Cartesian mesh $(ih, jh)$ $(i, j = 0, \ldots, n)$ of grid-spacing $h = 1/n$. We partition the index set $\{0, \ldots, n\}^2$ in the two disjoint sets $\Omega_D$ and $\Gamma_D$, where $\Gamma_D$ represents the discrete version of the set $\Gamma$. Finally, we use the scheme proposed in [RT92] to discretize $\|\nabla T\|$, arriving at the discrete equation

$$\max(D_{ij}^{-x}T, -D_{ij}^{+x}T, 0)^2 + \max(D_{ij}^{-y}T, -D_{ij}^{+y}T, 0)^2 = \frac{1}{F_{ij}^2}, \quad (i,j) \in \Omega_D \tag{3}$$
$$T_{ij} = 0, \quad (i,j) \in \Gamma_D.$$

Here $D_{ij}^{\pm x}T$ denotes the backward/forward finite-difference approximation of the partial derivative with respect to $x$, that is,

$$D_{ij}^{-x}T = \frac{T_{ij} - T_{i-1,j}}{h}, \quad -D_{ij}^{+x}T = \frac{T_{ij} - T_{i+1,j}}{h},$$

where $T_{ij} = T(ih, jh)$, $F_{ij} = F(ih, jh)$. Furthermore, we set $T_{i,n+1} = T_{i,-1} = \infty$ and $T_{n+1,j} = T_{-1,j} = \infty$ $(i, j = 0, \ldots, n)$ for notational convenience.

## 3. The Algorithm

The solution of (3) can be efficiently computed with the fast marching method (compare e.g. [Set99]), which is given as follows. First, all points in $\Gamma_D$ are tagged as *known*. All points in $\Omega_D$ that have a neighbor in $\Gamma_D$ are tagged as *trial*. Based on the known values in $\Gamma_D$, *trial values* are computed for all trial points, using (3). We denote by $\mathcal{N}$ the set of trial points, the so-called *narrow band*. Then the following iteration is performed:

(A1) Let $(i,j) \in \mathcal{N}$ denote the point with the smallest trial value.
(A2) Tag $(i,j)$ as *known*, and remove it from $\mathcal{N}$.
(A3) Tag all neighbors of $(i,j)$, which are not *known*, as *trial* (unless this hasn't already been done) and add them to $\mathcal{N}$.
(A4) Recompute the values of $T$ for all trial neighbors of $(i,j)$, using (3) with the *known* values of $T$ only.
(A5) If $\mathcal{N} \neq \emptyset$, goto (A1).

As the value of $T_{ij}$ in (3) depends only on smaller values of $T$ in the neighboring points, the fast marching method computes a solution to (3).

In the narrow band $\mathcal{N}$, the following estimate holds.

**Lemma 1.** *In every step of the fast marching method, there holds:*

$$\max_{(i,j) \in \mathcal{N}} T_{ij} - \min_{(i,j) \in \mathcal{N}} T_{ij} \leqslant \frac{1}{F_{\min}} \cdot h.$$

*Proof.* Let $\mathcal{N}$ denote the narrow band at the beginning of step (A1), and $\mathcal{N}^*$ the narrow band at the end of step (A4), after one cycle of the fast marching method. Let $(i,j) \in \mathcal{N}$ denote the point with the smallest trial value, and let $(k,l)$ denote some trial neighbor of $(i,j)$, with its trial value $T_{kl}$ computed in step (A4). Then

$$\frac{1}{F_{kl}^2} \geqslant \left(\frac{T_{kl} - T_{ij}}{h}\right)^2 \quad \Rightarrow \quad T_{kl} \leqslant \min_{(p,q) \in \mathcal{N}} T_{pq} + \frac{h}{F_{kl}}.$$

If we inductively assume the assertion to be true for $\mathcal{N}$, then, by the last inequality,

$$\max_{(p,q)\in\mathcal{N}^*} T_{pq} \leqslant \min_{(p,q)\in\mathcal{N}} T_{pq} + \frac{h}{F_{\min}}$$

(except for the neighbors of $(i,j)$, all grid-function values of points in $\mathcal{N}$ remain unchanged). On the other hand, the computation of a trial value, which depends on the new *known* value $T_{ij}$ cannot yield a smaller result than $T_{ij}$. Thus $\min_{(p,q)\in\mathcal{N}} T_{pq} \leqslant \min_{(p,q)\in\mathcal{N}^*} T_{pq}$, and we are done. □

The crucial point in the fast marching method is step $(A1)$. Typically, a heap based priority queue is used in order to keep track of the grid-point with the currently smallest trial value $T_{ij}$. Here we follow the suggestion of [YBS06] to use a so-called *untidy* priority queue within the fast marching method to reach linear run-time. The idea is to quantize the trial values contained in $\mathcal{N}$, and to pick any value, that quantizes to the minimal value, instead of the minimal value.

The untidy priority queue is organized as a circular array of $n_B + 1$ "buckets", denoted by $B_0, \ldots, B_{n_B}$, to store the narrow band $\mathcal{N}$. Those buckets are dynamical arrays, implemented as a singly-linked lists with FIFO ordering. Every bucket should contain *trial* points with *trial* values differing for at most $\delta$, where $\delta > 0$ is a small, appropriately chosen quantity. For the choice of the point with the smallest *trial* value in step (A1), we make no distinction between points contained in one bucket. In detail, the bucket $B_r$ contains all trial points $(i,j) \in \mathcal{N}$ with

$$r = \lfloor T_{ij}/\delta \rfloor \mod (n_B + 1), \quad \text{where } \delta = \frac{h}{F_{\min} \cdot n_B}. \qquad (4)$$

By the choice of $\delta$ it is ensured that only grid-points with trial values within the range of $\delta$ belong to a common bucket, that is $|T_{ij} - T_{kl}| < \delta$ for all $(i,j), (k,l) \in B_r$, as we have, similar to Lemma 1,

$$\max_{(i,j)\in\mathcal{N}} T_{ij} - \min_{(i,j)\in\mathcal{N}} T_{ij} < \frac{h}{F_{\min}} + \delta = \delta \cdot (n_B + 1)$$

$$\Rightarrow \left\lfloor \max_{(i,j)\in\mathcal{N}} T_{ij}\Big/\delta \right\rfloor - \left\lfloor \min_{(i,j)\in\mathcal{N}} T_{ij}\Big/\delta \right\rfloor < n_B + 1,$$

when using the untidy priority queue. On the other hand, the modulus operation in (4) has the effect that buckets, which are emptied during the algorithm, are refilled in subsequent cycles of the fast marching method for memory efficiency reasons. During the computation we have to keep track of the bucket $B_s$, which holds the grid-points with the smallest trial values. The untidy queue supports two operations, namely

- *Insertion of some point $(i,j)$:* To insert some point $(i,j)$, we compute $r$ by (4) and attach $(i,j)$ at the end of bucket $B_r$.
- *Deletion of the point $(i,j)$ with the approximately minimal trial value:* The index $s$ should hold the number of the bucket with the smallest trial value. If $B_s$ is empty, we search cyclically for the next non-empty bucket $B_{s'}$, passing cyclically from bucket to bucket, substituting $s$ by $\lfloor s+1 \rfloor \mod (n_B + 1)$. If all buckets are empty, the method terminates. Otherwise, we return the first element $(i,j)$ in $B_s$, and remove it from $B_s$.

Since the recomputation of some trial value $T_{kl}$ in step (A4) can yield a smaller value, it might become necessary to move $(k,l)$ to another bucket. In our approach we decided *to insert a grid-point once more* if it gets a new trial-value $\hat{T}_{kl}$ such that the corresponding bucket number changes. Additionally, if we encounter some already *known* point in the delete operation, we simply skip the point and pick the next one. This approach may have the effect that a grid-point appears up to four times in the untidy priority queue.[1] However, this way we avoid the costly removal of trial points from a bucket if they would have to be transferred to another one.

The complexity of the insert operation is obviously $\mathcal{O}(1)$. However, the computational cost of getting the next minimal element in the queue may not be constant in a single step, as it incorporates the search for a non-empty bucket. Nevertheless, we obtain linear complexity in total, as the following lemma shows.

**Lemma 2.** *The complexity of the modified fast marching method with $n_B$ buckets is $\mathcal{O}(N + n \cdot n_B)$, where $N = |\Omega_D|$.*

*Proof.* In every cycle (A1)-(A5) of the fast marching method, one grid-point becomes *known*, and after $N$ cycles the method terminates. Steps (A3) and (A4) require $\mathcal{O}(1)$ operations. But possibly several buckets have to be passed in the delete operation in step (A2). However, the solution of (3) is bounded from above by $M = 2/F_{\min}$. Let us imagine for simplicity, that we would use a linear array $B_0, B_1, B_2, \ldots$ of buckets for our untidy priority queue, disregarding the modulus operation in (4). Then we could cope with $\lfloor M/\delta \rfloor$ buckets. Hence also with the circular array, at most $M/\delta = 2/h \cdot n_B = \mathcal{O}(n \cdot n_B)$ buckets have to be traversed in the delete operation during the whole algorithm. $\square$

## 4. Estimates

Let us write in short $N_{ij}(T)$ for the approximation of $\|\nabla T(ih, jh)\|$ as given in formula (3). Then the following comparison principle holds.

**Lemma 3.** *Let $S_{ij}$ and $T_{ij}$ denote two grid-functions with $N_{ij}(S)F_{ij} \leqslant 1$ and $N_{ij}(T)F_{ij} \geqslant 1$ for all $(i,j) \in \Omega_D$, respectively. Then*

$$\max_{(i,j) \in \Omega_D} (S_{ij} - T_{ij})^+ \leqslant \max_{(i,j) \in \Gamma_D} (S_{ij} - T_{ij})^+$$

*Proof.* Assume that for some inner grid-point $(i,j) \in \Omega_D$ we have $S_{ij} - T_{ij} = \max_{(i,j) \in \Omega_D}(S_{ij} - T_{ij}) = \rho > 0$, and that additionally $(i,j)$ is the grid-point with the minimal value of $T_{ij}$ among all grid-points in $\Omega_D$ realizing that maximum. The solution $X$ of the discrete equation

$$\max(X - T_{i-1,j}, X - T_{i+1,j}, 0)^2 + \max(X - T_{i,j-1}, X - T_{i,j+1}, 0)^2 = \frac{h^2}{F_{ij}^2}. \quad (5)$$

at $(i,j)$ can be written as[2]

$$X = \min_{s+t=1,\, s,t \geqslant 0} \left\{ s \cdot \min(T_{i-1,j}, T_{i+1,j}) + t \cdot \min(T_{i,j-1}, T_{i,j+1}) + \frac{h}{F_{ij}}\sqrt{s^2 + t^2} \right\}.$$

Hence, by the assumption $N_{ij}(T)F_{ij} \geqslant 1$, we have $T_{ij} \geqslant X$, and there are numbers $s, t \geqslant 0$ with $s + t = 1$, and $\mu, \nu \in \{-1, 1\}$, such that

$$T_{ij} \geqslant s \cdot T_{i+\mu,j} + t \cdot T_{i,j+\nu} + \frac{h}{F_{ij}}\sqrt{s^2 + t^2}. \quad (6)$$

Similarly we deduce by $N_{ij}(S)F_{ij} \leqslant 1$ that

$$S_{ij} \leqslant s \cdot S_{i+\mu,j} + t \cdot S_{i,j+\nu} + \frac{h}{F_{ij}}\sqrt{s^2 + t^2}.$$

Taking the difference of the last two equations yields

$$\rho = S_{ij} - T_{ij} \leqslant s \cdot (S_{i+\mu,j} - T_{i+\mu,j}) + t \cdot (S_{i,j+\nu} - T_{i,j+\nu}) \leqslant \rho.$$

---

[1] In general, we observe that every point is inserted two times on average.

[2] In [BR06], a variational discretization of the eikonal equation directly leads to this formula. As easily verified by solving a scalar extreme value problem, $X$ fulfills (5).

By the choice of $\rho$, we would have $S_{i+\mu,j} - T_{i+\mu,j} = \rho$ if $s > 0$, and $S_{i,j+\nu} - T_{i,j+\nu} = \rho$ if $t > 0$. Together with (6) we deduce that $T_{ij}$ has a neighbor $T_{kl}$ with $T_{kl} < T_{ij}$ and $S_{kl} - T_{kl} = \rho$. By our assumption, $T_{kl}$ has to be a boundary point, thus the lemma is proved. □

With the help of the discrete comparison principle, we are now able to prove an estimate on the error caused by the inexact minimization in the untidy priority queue.

**Theorem 4.** *Let $\hat{T}$ denote the solution, computed with the modified fast marching method, and let $T$ denote the exact solution of (3). Then, for all $(i,j) \in \Omega_D$,*

$$0 \leqslant \frac{\hat{T}_{ij} - T_{ij}}{|\hat{T}_{ij}|} \leqslant \sqrt{2} \cdot \frac{F_{\max}}{h} \cdot \delta = \sqrt{2} \cdot \frac{F_{\max}}{F_{\min}} \cdot \frac{1}{n_B},$$

*where $n_B + 1$ denotes the number of buckets, and $\delta = \frac{h}{F_{\min} \cdot n_B}$ the increment.*

*Proof.* We analyze the error due to the untidy priority queue. When some value $\hat{T}_{ij}$ becomes *known*, then there might be a smaller value $\hat{T}_{kl}$ in the same bucket that would have been accepted before $\hat{T}_{ij}$ if an exact priority queue had been used, with $\hat{T}_{kl} \leqslant \hat{T}_{ij} < \hat{T}_{kl} + \delta$. Let $X$ denote the solution of

$$\max(X - \hat{T}_{i-1,j}, X - \hat{T}_{i+1,j}, 0)^2 + \max(X - \hat{T}_{i,j-1}, X - \hat{T}_{i,j+1}, 0)^2 = \frac{h^2}{F_{ij}^2}.$$

Of course, we have $X \leqslant \hat{T}_{ij}$. Since, as we have seen, all values $\hat{T}_{kl}$ that are accepted after $\hat{T}_{ij}$ fulfill $\hat{T}_{kl} > \hat{T}_{ij} - \delta$, we also have $X > \hat{T}_{ij} - \delta$. Hence,

$$\frac{h^2}{F_{ij}^2} \leqslant \max(\hat{T}_{ij} - \hat{T}_{i-1,j}, \hat{T}_{ij} - \hat{T}_{i+1,j}, 0)^2 + \max(\hat{T}_{ij} - \hat{T}_{i,j-1}, \hat{T}_{ij} - \hat{T}_{i,j+1}, 0)^2$$

$$< \left[\max(X - \hat{T}_{i-1,j}, X - \hat{T}_{i+1,j}, 0) + \delta\right]^2 + \left[\max(X - \hat{T}_{i,j-1}, X - \hat{T}_{i,j+1}, 0) + \delta\right]^2$$

$$\leqslant \left(\frac{h}{F_{ij}} + \sqrt{2} \cdot \delta\right)^2 \leqslant \frac{h^2}{F_{ij}^2}\left(1 + \sqrt{2} \cdot \delta \cdot \frac{F_{\max}}{h}\right)^2.$$

Thus $N_{ij}(\hat{T})F_{ij} \geqslant 1$, and as $N_{ij}(T)F_{ij} = 1$, we obtain from Lemma 3 that $T_{ij} \leqslant \hat{T}_{ij}$ for all $(i,j)$. On the other hand, we have $\theta \cdot N_{ij}(\hat{T})F_{ij} = N_{ij}(\theta \cdot \hat{T})F_{ij} \leqslant 1$ with $\theta = (1 + \sqrt{2}\delta F_{\max}/h)^{-1}$, $0 < \theta < 1$. Once again, Lemma 3 yields

$$0 \leqslant \hat{T}_{ij} - T_{ij} \leqslant (1 - \theta) \cdot \hat{T}_{ij} + \max_{(k,l) \in \Omega_D}(\theta \hat{T}_{kl} - T_{kl})^+ \leqslant \sqrt{2} \cdot \frac{F_{\max}}{h} \cdot \delta \cdot \hat{T}_{ij}.$$

□

**Remarks:**

(1) Of course, the above result can be generalized to arbitrary space dimensions. In the estimate on the relative error, the factor $\sqrt{2}$ has then to be substituted by the square root of the space dimension.
(2) As Lemma 2 and Theorem 4 show, an increase in the number of buckets causes a lower error at the price of a higher complexity. Moreover, we expect a total error of at least $\mathcal{O}(h)$ due to the discretization of the eikonal equation, and a complexity of at least $\mathcal{O}(N)$, as $N$ is the number of gridpoints. Hence, the suitable choice for $n_B$ is $n_B = \mathcal{O}(n) = \mathcal{O}(1/h)$. A choice of $n_B = F_{\max}/F_{\min} \cdot n$ buckets would yield an error bound independent of the speed ratio $F_{\max}/F_{\min}$, but a total complexity of $\mathcal{O}(F_{\max}/F_{\min} \cdot N)$.

(3) The untidy priority queue can be used to modify the fast marching method on triangulations [KS98], the method for the solution of the eikonal equation on parametric manifolds [SK04], or the ordered upwind method for more general Hamilton-Jacobi equations given in [SV01], to reduce the run-time to $\mathcal{O}(N)$.

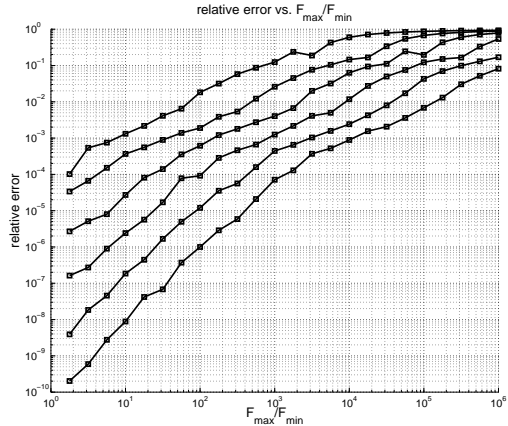

FIGURE 1. Relative error vs. the ratio $F_{max}/F_{min}$ for different numbers of buckets (from top to bottom we used 2, 8, 32, 128, 512 and 2048 buckets).

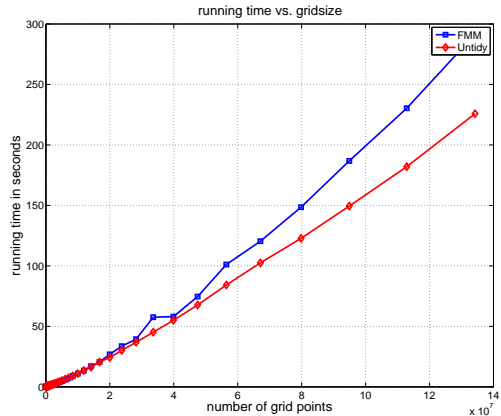

FIGURE 2. Run-time compared to the original fast marching method.

## 5. EXAMPLES

The first example concerns the error estimation asserted in Theorem 4. We used the speed function $F(x) = (1 + \frac{r-1}{2}) + \frac{r-1}{2}\sin(2\pi |x|)$ with $r = F_{max}/F_{min}$ on $\Omega = [0,1]^2$ and the boundary value $T(1,0) = 0$. Figure 1 illustrates the linear dependence of the relative error $(\hat{T}_{ij} - T_{ij})/\hat{T}_{ij}$ on the ratio $F_{max}/F_{min}$ for different numbers of buckets.

The second example compares the run-time of the algorithm with an untidy priority queue to the original fast marching Method. We computed $T_{ij}$ on $\Omega = [0,1]^2$ with the speed function $1/F$ sampled from a $[0,1]-$uniform distributed random variable. Figure 2 shows an asymptotically linear run-time of the modified method, and a non-linear behavior of the original method.

*Acknowledgment:* We would like to thank Folkmar Bornemann for his helpful comments on our paper.

## References

[BR06]  Folkmar Bornemann and Christian Rasch, *Finite-element discretization of static Hamilton-Jacobi equations based on a local variational principle.*, Comput. Visual Sci. **9** (2006), 57–69. 4

[KS98]  Ron Kimmel and James A. Sethian, *Computing geodesic paths on manifolds.*, Proc. Natl. Acad. Sci. USA **95** (1998), no. 15, 8431–8435. 6

[RT92]  Elisabeth Rouy and Agnés Tourin, *A viscosity solutions approach to shape-from-shading.*, SIAM J. Numer. Anal. **29** (1992), no. 3, 867–884. 2

[Set96]  James A. Sethian, *A fast marching level set method for monotonically advancing fronts.*, Proc. Natl. Acad. Sci. USA **93** (1996), no. 4, 1591–1595. 1

[Set99]  ———, *Fast marching methods.*, SIAM Rev. **41** (1999), no. 2, 199–235. 2

[SK04]  Alon Spira and Ron Kimmel, *An efficient solution to the eikonal equation on parametric manifolds.*, Interfaces Free Bound. **6** (2004), no. 3, 315–327. 6

[SV01]  James A. Sethian and Alexander Vladimirsky, *Ordered upwind methods for static Hamilton-Jacobi equations.*, Proc. Natl. Acad. Sci. USA **98** (2001), no. 20, 11069–11074. 6

[Tsi95]  John N. Tsitsiklis, *Efficient algorithms for globally optimal trajectories.*, IEEE Trans. Autom. Control **40** (1995), no. 9, 1528–1538. 1

[YBS06]  Liron Yatziv, Alberto Bartesaghi, and Guillermo Sapiro, *O(N) implementation of the fast marching algorithm.*, J. Comput. Phys. **212** (2006), no. 2, 393–399. 1, 3

*E-mail address*: `rasch@ma.tum.de, satzger@in.tum.de`

Center of Mathematics, Technical University of Munich, 80290 Munich, Germany